\documentclass[12pt]{article}
\usepackage{amsfonts,amssymb,amsmath,theorem}

\newtheorem{theorem}{Theorem}
\newtheorem{lemma}[theorem]{Lemma}
\newtheorem{prop}[theorem]{Proposition}

\newcommand\proof[1]{\noindent\textit{Proof#1}}
\newcommand\proofskip{\vspace{\theorempostskipamount}}

\let\eps=\varepsilon
\newcommand\Int{\mathop{\mathrm{Int}}}
\renewcommand\circ{\mspace{2mu}{\mathchar"220E}\mspace{2mu}}

\newcommand\abs[1]{\left|#1\right|}              
\newcommand\norm[1]{\left\|#1\right\|}           
\newcommand\qed{\ifhmode\unskip\nobreak\fi\quad  
   \ifmmode\square\else\hbox{$\square$}\fi}      

\begin{document}

\begin{center}
\large\bfseries
T-ENTROPY AND VARIATIONAL PRINCIPLE\\
FOR THE SPECTRAL RADIUS\\
OF WEIGHTED SHIFT OPERATORS

\bigskip\medskip\normalsize\mdseries
V.\,I.\ BAKHTIN

\medskip\itshape
Belarus State University (e-mail: bakhtin@tut.by)

\end{center}

\renewcommand\abstractname{}
\begin{abstract}
In this paper we introduce a new functional invariant of discrete time dynamical systems --- the
so-called \emph{$t$-entropy}. The main result is that this $t$-entropy is the Legendre dual
functional to the logarithm of the spectral radius of the weighted shift operator on $L^1(X,m)$
generated by the dynamical system. This result is called the \emph{Variational principle} and is
similar to the classical variational principle for the topological pressure.
\end{abstract}

\bigbreak\bigskip

\quad\parbox{14.0cm} {\textbf{Keywords:} {\itshape weighted shift operator, variational principle,
$t$-entropy, entropy\\ statistic theorem}

\medbreak \textbf{2000 MSC:} 37A35, 47B37 }

\bigbreak\medskip\bigskip

In the present paper we define and investigate a new functional invariant of discrete time
dynamical systems. It is called \emph{$t$-entropy}. Similarly to the classical Kolmogorov--Sinai
entropy this invariant is a concave functional on the set of probability measures defined on the
phase space of a dynamical system. It is well known that the Fenchel--Legendre transform of the
Kolmogorov--Sinai entropy coincides with the topological pressure (this fact is usually called the
`variational principle for the topological pressure'). Similarly,  the basic property of
$t$-entropy is that its Fenchel--Legendre transform turns out to be equal to the logarithm of the
spectral radius of a weighted shift operator generated by the dynamical system. Moreover,
$t$-entropy is always upper semicontinuous whereas the Kolmogorov--Sinai entropy in general does
not possess this property.

The problem of existence of an invariant that would make it possible to prove the variational
principle for the spectral radius of a weighted shift operator was set up by A.\,V.\ Lebedev. He
also proposed the term `$t$-entropy' itself. This problem was under discussion at the seminar
directed by A.\,B.\ Antonevich and A.\,V.\ Lebedev in the course of several years and finally was
solved only due to their constant inspiring influence. Initial Lebedev's idea was very simple: to
proceed in the same manner as R.\ Bowen in \cite{1} where he proved the variational principle for
the topological pressure. But the further analysis showed that Bowen's method does not work in the
new setting. So it was radically revised.

In essence, the main result of this paper was announced in \cite{2, 3, 4} (in a slightly different
form and in a bit more general setting). Nevertheless the proof of it in the present setting  is
published now for the first time.

The paper is organized as follows. In Section \ref{1..} we define the notion of $t$-entropy,
formulate the Variational principle, and prove its easy part. In Section \ref{2..} we study basic
properties of $t$-entropy. In Section~\ref{3..} we formulate the so-called Entropy Statistic
Theorem (Theorem \ref{..2}) and deduce from it the hard part of the Variational principle. Section
\ref{4..} contains the proof of the Entropy Statistic Theorem.

\section{Variational Principle}\label{1..}

Let  $\alpha\!: X\to X$ be  a measurable mapping of a measurable space $(X,\mathfrak A)$ supplied
with a $\sigma$-finite measure $m$. This mapping generates the shift operator $A$ that maps each
function $f$ to $Af =f\circ\alpha$. Assume that $A$ is bounded on the space $L^1(X,m)$ of
integrable real-valued functions on $X$. The boundedness condition for $A$ is equivalent to the
existence of a constant $C$ such that $m(\alpha^{-1}(G))\le Cm(G)$ for every  $G\in \mathfrak A$.
The latter inequality implies that  $\norm A\le C$. For each function $\varphi\in L^\infty(X,m)$,
where $L^\infty(X,m)$ is the space of essentially bounded real-valued functions, we define the
weighted shift operator $A_\varphi$ acting on $L^1(X,m)$ by the formula
\begin{equation}\label{,,1}
\bigl[A_\varphi f\bigr](x) =e^{\varphi(x)}f(\alpha(x)),\qquad f\in L^1(X,m).
\end{equation}
Denote by $\lambda(\varphi)$ the logarithm of the spectral radius of $A_\varphi$:
\begin{equation}\label{,,2}
\lambda(\varphi) =\lim_{n\to\infty}\frac{1}{n}\ln\norm{A_\varphi^n}.
\end{equation}

A linear functional on $L^\infty(X,m)$ is called {\em positive\/} if it is nonnegative on the set
of all nonnegative functions, and {\em normalized\/} if it takes the unit value at the unit
function. Denote by $M(X,m)$ the set of all positive normalized linear functionals on
$L^\infty(X,m)$. Obviously, $M(X,m)$ can be identified with the set of all finitely additive
probability measures on $\mathfrak A$ that are absolutely continuous with respect to $m$. Therefore
we will call the elements of $M(X,m)$ {\em measures}. Any finite set $D =\{g_1,\dots,g_k\}$ of
nonnegative functions $g_i\in L^\infty(X,m)$ will be called a {\em measurable partition of unity on
$X$\,} if it satisfies (almost everywhere) the identity $g_1+ \dots +g_k\equiv 1$. For any measure
$\mu\in M(X,m)$ we define the {\em $t$-entropy\/} $\tau(\mu)$ by the formulae
\begin{gather}\label{,,3}
\tau(\mu) =\inf_{n\in\mathbb N}\frac{\tau_n(\mu)}{n},\qquad
\tau_n(\mu) =\inf_D\tau_n(\mu,D),\\[6pt]
\tau_n(\mu,D) =\sup_{\norm f =1}\sum_{g\in D}\mu(g) \ln\frac{\int_X
g\abs{f\circ\alpha^n}dm}{\mu(g)}. \label{,,4}
\end{gather}
Infimum in \eqref{,,3} is taken over all measurable partitions of unity $D$ on $X$ and the supremum
in \eqref{,,4} is taken over all $f\in L^1(X,m)$. If $\mu(g) =0$ for some $g\in D$, then the
corresponding summand in~\eqref{,,4} is assumed  to be zero regardless the integral in numerator.
But if there exists a function $g\in D$ such that $\int_X g\,dm =0$ and $\mu(g) >0$, then we set
$\tau(\mu) =-\infty$.

A measure $\mu\in M(X,m)$ is called \emph{$\alpha$-invariant} if $\mu(f\circ\alpha) =\mu(f)$ for
all $f\in L^\infty(X,m)$. This is equivalent to the condition $\mu(\alpha^{-1}(G)) =\mu(G)$,\ \
$G\in\mathfrak A$. We denote by $M_\alpha(X,m)$ the set of all \hbox{$\alpha$-}invariant measures
$\mu \in M(X,m)$.

\begin{theorem}[Variational principle]\label{..1}
The logarithm of the spectral radius of the weighted shift operator\/ \eqref{,,1} satisfies the
variational principle
\begin{equation}\label{,,5}
\lambda(\varphi) =\max_{\mu\in M_\alpha(X,m)} \bigl(\mu(\varphi) +\tau(\mu)\bigr), \qquad
\varphi\in L^\infty(X,m).
\end{equation}
\end{theorem}

This is the main result of the paper.  Let us make at once some relevant remarks. First,
formula~\eqref{,,5} means that $\lambda(\varphi)$ is the Fenchel--Legendre transform of the
restriction of the functional $-\tau(\mu)$ to the domain $M_\alpha(X,m)$. Secondly, if the maximum
in \eqref{,,5} is attained at a measure $\mu$ then $\mu$ is a subgradient of~$\lambda(\varphi)$.
This implies that $\lambda(\varphi)$ is convex with respect to  $\varphi$.

\section{Proof of the Variational Principle}\label{2..}

First let us verify the inequality $\lambda(\varphi) \ge \mu(\varphi) +\tau(\mu)$, where $\mu\in
M_\alpha(X,m)$. To this end introduce the notation $S_n\varphi =\varphi +\varphi\circ\alpha+
\,\dotsm\,+ \varphi\circ \alpha^{n-1}$. Then \eqref{,,1} implies the equality $A_\varphi^n f
=e^{S_n\varphi} f\circ \alpha^n$. For an arbitrary $n\in\mathbb N$ and $\eps>0$ choose a measurable
partition of unity $D$ such that for any function $g\in D$ the essential oscillation of
$S_n\varphi$ over the support of $g$ does not exceed $\eps$. Set $S_n\varphi(g) =
\mathop{\mathrm{ess\,sup}}\{\,S_n\varphi(x)\mid g(x)\ne 0\,\}$. Then for any $\mu\in M_\alpha(X,m)$
and $f\in L^1(X,m)$ by the concavity of the logarithm function we have
\begin{align*}
\eps+\ln\int_X e^{S_n\varphi}\abs{f\circ\alpha^n}dm &\ge \ln\sum_{g\in D}e^{S_n\varphi(g)}
\int_X g\abs{f\circ\alpha^n}dm\\[3pt]
&\ge \sum_{\mu(g)\ne 0}\mu(g)\ln\frac{e^{S_n\varphi(g)}
\int_X g\abs{f\circ\alpha^n}dm}{\mu(g)}\\[3pt]
&\ge \mu(S_n\varphi) +\sum_{\mu(g)\ne 0}\mu(g)\ln \frac{\int_X g\abs{f\circ\alpha^n}dm}{\mu(g)}.
\end{align*}
Pass here to the supremums over all $f$ with unit norm. As a result, taking into
account~\eqref{,,4}, we obtain the inequality $\eps+ \ln\norm{A_\varphi^n}\ge \mu(S_n\varphi)
+\tau_n(\mu,D)$. Divide it by $n$ and turn $n$ to infinity. The limit will be $\lambda(\varphi)\ge
\mu(\varphi) +\tau(\mu)$.

It is much more difficult  to prove that $\lambda(\varphi)$ does not exceed the right-hand side
of~\eqref{,,5}. We will deduce this fact from the Entropy Statistic Theorem for empirical measures,
which will be formulated below.

Take any point $x\in X$. The empirical measures $\delta_{x,n}\in M(X,m)$ are defined by the formula
\begin{equation*}
\delta_{x,n}(f) =\frac{1}{n}S_nf(x) = \frac{1}{n}\Bigl(f(x) +f(\alpha(x))+\,\dotsm\,
+f\bigl(\alpha^{n-1}(x)\bigr)\Bigr).
\end{equation*}

\smallskip\noindent
Since here $f\in L^\infty(X,m)$, the value of $\delta_{x,n}(f)$ cannot be uniquely determined for
each~$x$. Indeed, this value, as a function of~$x$, is an element of $L^\infty (X,m)$. So it is
defined only up to a set of measure zero.

Supply $M(X,m)$ with a *-weak topology of the dual space to $L^\infty (X,m)$. Given a measure
$\mu\in M(X,m)$ and its neighborhood $O(\mu)$ let us consider the sequence of sets
\begin{equation*}
X_n(O(\mu)) =\{\,x\in X\mid \delta_{x,n}\in O(\mu)\}.
\end{equation*}
Of course, they are defined only up to a set of measure zero. Nevertheless, we can state

\begin{theorem}[Entropy Statistic Theorem]\label{..2}
For any\/ $\mu\in M(X,m)$ and any\/ $t>\tau(\mu)$ there exist a neighborhood\/~$O(\mu)$ in the\/
*-weak topology and a large enough number\/ $C(t,\mu)$ such that for all\/ $f\in L^1(X,m)$ and\/
$n\in \mathbb N$ we have
\begin{equation}\label{,,6}
\int_{X_n(O(\mu))} f\circ\alpha^n\,dm\le C(t,\mu)e^{nt}\int_X\abs{f}\,dm, \qquad t>\tau(\mu).
\end{equation}
\end{theorem}

This very hard theorem will be proved in Section \ref{4..}. Now we are using it to complete the
proof of the Variational principle. We start with two auxiliary propositions.

\begin{prop}\label{..3}
The sets\/ $M(X,m)$ and\/ $M_\alpha(X,m)$ are compact in the\/
*-weak topology.
\end{prop}

\proof. These are closed subsets of the unit ball in the dual space  to $L^\infty(X,m)$. This ball
is compact by Alaoglu's theorem. \qed

\begin{prop}\label{..4}
If\/ $U$ is a neighborhood of\/ $M_\alpha(X,m)$ in\/ $M(X,m)$ then there exists so large\/ $N$
that\/ $\delta_{x,n} \in U$ for almost all\/ $x\in X$ and all\/ $n>N$.
\end{prop}

\proof. By the previous Proposition the set $M(X,m)\setminus U$ is compact. Hence there exist a
finite collection of functions $f_1,\,\dotsc,\,f_k$ from $L^\infty(X,m)$ and a positive number
$\eps$ such that $\sum_{i=1}^k\bigl|\mu(f_i) -\mu(f_i\circ\alpha) \bigr| >\eps$ for all $\mu\in
M(X,m)\setminus U$. Obviously, $\delta_{x,n}(f_i) -\delta_{x,n}(f_i\circ\alpha) = n^{-1}\bigl(f(x)
-f(\alpha^{n}(x))\bigr)$. Therefore for all $n$ large enough we have (almost everywhere)
\begin{equation*}
\sum_{i=1}^k\bigl|\delta_{x,n}(f_i) -\delta_{x,n} (f_i\circ\alpha) \bigr| \le
\frac{2}{n}\sum_{i=1}^k \mathop{\mathrm{ess\,sup}}\abs{f_i} <\eps.
\end{equation*}
From here it follows that $\delta_{x,n}\in U$. \qed \proofskip

Now we are able to finish the proof of the Variational principle.

\proofskip\noindent
Set
\begin{equation}\label{,,7}
\Lambda(\varphi) =\sup_{\mu\in M_\alpha(X,m)} \bigl(\mu(\varphi) +\tau(\mu)\bigr).
\end{equation}

\medskip\noindent
Let us fix arbitrary  numbers $c>\Lambda(\varphi)$ and $\eps>0$ and define the functional $t(\mu) =
c-\mu(\varphi)$ on the set of invariant measures $M_\alpha(X,m)$. Clearly, $t(\mu)> \tau(\mu)$. For
each $\mu\in M_\alpha(X,m)$ choose a neighborhood~$O(\mu)$ in $M(X,m)$ such that the Entropy
Statistic Theorem holds true for this $O(\mu)$ and $t=t(\mu)$ and at the same time for all $\nu\in
O(\mu)$ we have the estimate $\nu(\varphi) <\mu(\varphi) +\eps$. Then for almost all $x\in
X_n(O(\mu))$ the next relations hold
\begin{equation}\label{,,8}
S_n\varphi(x) =n\delta_{x,n}(\varphi) <n\bigl(\mu( \varphi) +\eps\bigr) = n\bigl(c-t(\mu)
+\eps\bigr).
\end{equation}

By Proposition \ref{..3} the set $M_\alpha(X,m)$ is compact. Let us cover it by a finite collection
of neighborhoods $O(\mu_1),\,\dots,\, O(\mu_k)$ of the form described above. It follows from
Proposition~\ref{..4} that for all $n$ large enough almost every point $x\in X$ belongs to at least
one of the sets $X_n(O(\mu_i))$,\ \ $i=1,\,\dotsc,\,k$. Therefore, using \eqref{,,8} and the
Entropy Statistic Theorem we obtain the following estimate for each function $f\in L^1(X,m)$ with
unit norm
\begin{align*}
\norm{A_\varphi^n f} &=\int_X\abs{e^{S_n\varphi}f\circ\alpha^n}dm \le
\sum_{i=1}^k\int_{X_n(O(\mu_i))}e^{S_n\varphi}\abs{f}\circ
\alpha^n\,dm\\[3pt]
&\le \sum_{i=1}^k e^{n(c-t(\mu_i)+\eps)} C(t(\mu_i),\mu_i)e^{nt(\mu_i)} = e^{n(c+\eps)}
\sum_{i=1}^k C(t(\mu_i),\mu_i).
\end{align*}
When $n\to\infty$ this implies $\lambda(\varphi)\le c+\eps$. And since the numbers
$c>\Lambda(\varphi)$ and $\eps>0$ are arbitrary it follows  that $\lambda(\varphi) \le
\Lambda(\varphi)$. Finally, in Proposition \ref{..7} in the next section it will be proved that the
$t$-entropy $\tau(\mu)$ is upper semicontinuous on the compact set $M_\alpha(X,m)$. So the supremum
in \eqref{,,7} is, in fact, maximum. This completes the proof of Theorem \ref{..1}. \qed

\section{Properties of  $\mathbf t$-entropy}\label{3..}

In this section we study some properties of  $t$-entropy and in particular its  upper
semicontinuity. First note that $\tau(\mu)\le \lambda(0)$. Indeed, from \eqref{,,4} and the
concavity of the logarithm function it follows that for any measurable partition of unity $D$ on
$X$ we have
\begin{align*}
\tau_n(\mu,D) &=\sup_{\norm f =1}\sum_{g\in D}\mu(g) \ln\frac{\int_X g\abs{f\circ\alpha^n}dm}
{\mu(g)}\le \sup_{\norm f=1}\ln\sum_{g\in D}\int_X g\abs{f\circ\alpha^n}dm \\[3pt]
&= \sup_{\norm f=1}\ln \int_X \abs{f\circ\alpha^n}dm =\ln\norm{A^n}
\end{align*}
and hence, taking into account \eqref{,,3},
\begin{equation*}
\tau(\mu) =\inf_{n,D}\frac{\tau_n(\mu,D)}{n} \le \inf_n \frac{\ln\norm{A^n}}{n} =\lambda(0).
\end{equation*}
In particular, if the measure $m$ is invariant then $\norm{A^n} =1$ and therefore $\tau(\mu)\le 0$.

\begin{prop}\label{..5}
The functions\/ $\tau_n(\mu,D)$ and\/ $\tau(\mu)$ are concave with respect to\/ $\mu\in M(X,m)$.
\end{prop}

\proof. Suppose $\mu_1,\mu_2\in M(X,m)$ and $\mu =p_1\mu_1+p_2\mu_2$, where $p_1+p_2 =1$ and
$p_1,p_2\ge 0$. Then
\begin{gather*}
p_1\mu_1(g)\ln\frac{\int_X g\abs{f_1\circ\alpha^n}dm}{\mu_1(g)} +p_2\mu_2(g)\ln\frac{\int_X
g\abs{f_2\circ\alpha^n}dm}{\mu_2(g)}
\\[6pt]
\le \mu(g)\ln\frac{\int_X g\,(p_1\abs{f_1}+p_2\abs{f_2})\circ \alpha^n\,dm}{\mu(g)}.
\end{gather*}

\smallskip\noindent
Let us sum this inequality over $g\in D$ and pass to the supremums over $f_1$ and~$f_2$. As a
result we obtain the inequality $p_1\tau_n(\mu_1,D) +p_2\tau_n(\mu_2,D)\le \tau_n(\mu,D)$. It means
exactly that $\tau_n(\mu,D)$ is concave with respect to $\mu$. So \eqref{,,3} implies that
$\tau(\mu)$ is concave as well. \qed \proofskip

Let  $D =\{g_1,\dots,g_k\}$ be  a measurable partition of unity on $X$. Let us remove from it all
the elements~$g_i$ such that $\int_X g_i\,dm =0$ and put $D_m = \bigl\{g\in D\bigm| \int_X g\,dm
>0\bigr\}$. Denote by $M(D)$ and $M(D_m)$ the sets of all probability measures on the finite sets
$D$ and~$D_m$, respectively. Obviously, the sets $M(D)$ and~$M(D_m)$ are finite-dimensional
simplexes. If one extends each measure $\mu\in M(D_m)$ to $D\setminus D_m$ by zero, then the
simplex $M(D_m)$ becomes a face of $M(D)$. Note that formula \eqref{,,4} defines the functions
$\tau_n(\mu,D)$ not only for the measures $\mu\in M(X,m)$ but for $\mu\in M(D)$ as well.

\begin{prop}\label{..6}
The function\/ $\tau_n(\,\cdot\,,D)$ is continuous on\/ $M(D_m)$ and turns into\/~$-\infty$ on\/
$M(D)\setminus M(D_m)$.
\end{prop}

\proof. Let us fix a nonnegative function $h\in L^1(X,m)$ with unit norm
 and such that $\int_X g\,h\circ\alpha^n\,dm >0$ for all $g\in D_m$.
For any nonnegative function $f\in L^1(X,m)$ with unit norm let us define the family $f_\eps
=(1-\eps)f +\eps h$ depending on the parameter $\eps\in [0,1]$. In addition we introduce the
notation
\begin{equation*}
\psi(\mu,f) =\sum_{g\in D_m} \mu(g)\ln\int_X g\,f\circ\alpha^n\,dm, \qquad \psi_\eps(\mu)
=\sup_{\norm f=1}\psi(\mu,f_\eps).
\end{equation*}
Then for any $\mu,\nu\in M(D_m)$ we have
\begin{equation}\label{,,9}
\bigl|\psi(\mu,f_\eps) -\psi(\nu,f_\eps)\bigr|\le \max_{g\in D_m} \biggl|\ln\int_X
g\,f_\eps\circ\alpha^n\,dm\biggr| \sum_{g\in D_m}\bigl|\mu(g) -\nu(g)\bigr|.
\end{equation}
By construction,
\begin{equation}\label{,,10}
\int_X g\,f_\eps\circ\alpha^n\,dm \ge \eps\int_X g\,h\circ\alpha^n\,dm >0
\end{equation}
and on the other hand
\begin{equation}\label{,,11}
\int_X g\,f_\eps \circ\alpha^n\,dm\le \norm{A^n}.
\end{equation}

\medskip\noindent
From \eqref{,,9}--\eqref{,,11} it follows that for strictly positive $\eps$ the function
$\psi_\eps(\mu)$ depends continuously on the measure $\mu\in M(D_m)$. Evidently,
\begin{equation*}
\psi_0(\mu)\ge \psi_\eps(\mu)\ge \psi_0(\mu)+\ln (1-\eps).
\end{equation*}
Hence the function $\psi_0(\mu)$ is the uniform limit (as $\eps\to 0$) of $\psi_\eps(\mu)$ and so
it is also continuous. Furthermore, the difference of two continuous functions $\psi_0(\mu)\,
-\,\sum_{g\in D} \mu(g)\ln \mu(g)$ coincides with $\tau_n(\mu,D)$ and is continuous on
the~$M(D_m)$. The second part of Proposition~\ref{..6} follows immediately from~\eqref{,,4}. \qed
\proofskip

By virtue of Proposition \ref{..6} the function $\tau_n(\mu,D)$ is upper semicontinuous (in the
*-weak topology) on the $M(X,m)$. The same is true for the $t$-entropy $\tau(\mu) =
\inf_{n,D}\tau_n(\mu,D)/n$. So we have obtained

\begin{prop}\label{..7}
The functional\/ $\tau(\mu)$ is upper semicontinuous on the\/~$M(X,m)$.
\end{prop}

In the next four Propositions we observe some additional properties of the $t$-entropy although
they are not used in verification of the Variational principle.

For any partition of unity $D$ on $X$ and a positive number $\eps$ denote by $W(D,\eps)$ the set of
all measurable partitions of unity $E$ that satisfy the following condition: for each $h\in E$ the
oscillation of any function $g\in D$ on the support of $h$ does not exceed $\eps$.

\begin{prop}\label{..8}
Given a measure\/ $\mu\in M(X,m)$, a measurable partition of unity\/ $D$ on\/ $X$, and a number\/
$t>\tau_n(\mu,D)$ there exists a small\/ $\eps>0$ such that for all\/ $E\in W(D,\eps)$ the
estimate\/ $\tau_n(\mu,E) <t$ is true.
\end{prop}

\proof. Denote by $D_\mu$ the collection of all functions $g\in D$ satisfying the condition
$\mu(g)>0$. Choose a positive $\eps$ so small that for all $g_0\in D_\mu$,
\begin{equation}\label{,,12}
\mu(g_0)\ln\frac{\sqrt\eps +\eps\norm{A}^n}{\mu(g_0)}\, + \sum_{g\in D_\mu\setminus\{g_0\}}
\mu(g)\ln \frac{2\norm{A}^n}{\mu(g)} <t.
\end{equation}
Consider any partition of unity $E\in W(D,\eps)$. Let $E_\mu = \{\,h\in E\mid \mu(h)>0\,\}$. Then
for each nonnegative function $f\in L^1(X,m)$ with unit norm,  the concavity of the logarithm
implies the following relations
\begin{gather}\notag
\sum_{h\in E_\mu}\mu(h)\ln\frac{\int_X h\,f\circ\alpha^n\,dm}{\mu(h)} =\sum_{g\in
D_\mu}\mu(g)\sum_{h\in E_\mu}\frac{\mu(gh)}{\mu(g)}\ln
\frac{\int_X h\,f\circ\alpha^n\,dm}{\mu(h)} \displaybreak[0] \\[3pt]\notag
\le \sum_{g\in D_\mu}\mu(g)\ln\frac{1}{\mu(g)}\sum_{h\in E_\mu}
\frac{\mu(gh)}{\mu(h)}\int_X h\,f\circ\alpha^n\,dm \\[3pt]
\le \sum_{g\in D_\mu}\mu(g)\ln\frac{\int_X (g+\eps)f\circ \alpha^n\,dm}{\mu(g)}. \label{,,13}
\end{gather}
To finish the proof it  suffice  to check that \eqref{,,13} does not exceed $t$ provided $\eps$ is
small. Since, once it is proved then passing to the supremum over $f$ we obtain immediately the
desired estimate $\tau_n(\mu,E)\le t$.

If there exists a function $g_0\in D_\mu$ such that $\int_X g_0f\circ\alpha^n\,dm <\sqrt\eps$ then,
obviously, $\int_X (g_0+\eps) f\circ\alpha^n\,dm <\sqrt\eps+\eps\norm{A}^n$ and so, in view of
\eqref{,,12}, the whole of sum~\eqref{,,13} does not exceed $t$. In the opposite  case for each
$g\in D_\mu$ we have the estimate
\begin{equation*}
\int_X (g+\eps)f\circ\alpha^n\,dm\le \bigl(1+\sqrt\eps\norm{A}^n\bigr) \int_X
g\,f\circ\alpha^n\,dm.
\end{equation*}
Therefore, it follows  that \eqref{,,13} does not exceed the sum $\tau_n(\mu,D)
+\ln\bigl(1+\sqrt\eps\norm{A}^n\bigr)$, which is less then $t$ provided $\eps$ is small enough.
\qed

\begin{prop}\label{..9}
If one uses in\/ \eqref{,,3} only those partitions of unity\/ $D$ that consist of index functions
(corresponding to finite measurable partitions of the space\/ $X$) then the functionals\/
$\tau_n(\mu)$ and $\tau(\mu)$ do not change.
\end{prop}

\proof. This follows from Proposition \ref{..8}. \qed

\begin{prop}\label{..10}
If the mapping\/ $\alpha$ is invertible and\/ $\alpha^{-1}$ is measurable and the measure\/ $m$
is\/ $\alpha$-invariant then\/ $\tau(\mu)\equiv 0$.
\end{prop}

\proof. As it was already proved, if $m$ is invariant then $\tau(\mu) \le 0$. Consider any
measurable partition of unity $D$ on $X$ consisting of index functions of a finite family of  sets
$G_1,\,\dots,\, G_k$ which is  a measurable partition of $X$. One can assign to each~$G_i$ a
nonnegative measurable function $f_i$ such that it vanishes outside $G_i$ and at the same time
$\int_{G_i} f_i\,dm =\mu(G_i)$. Consider the function $f=(f_1+\dots+f_k)\circ \alpha^{-n}$. It is
easily seen that $\norm f =1$. Substituting~$f$ in \eqref{,,4} we obtain the inequality
$\tau_n(\mu,D) \ge 0$. Then by the previous proposition $\tau_n(\mu) \ge 0$ and so~$\tau(\mu) \ge
0$. \qed

\begin{prop}\label{..11}
If\/ $\mu\in M_\alpha(X,m)$ then\/ $\tau_{n+k} (\mu)\le \tau_n(\mu) +\tau_k(\mu)$. So for an\/
$\alpha$-invariant measure\/ $\mu$  \,$t$-entropy may be defined as the limit\/ $\tau(\mu)
=\lim_{n\to\infty}\tau_n(\mu)/n$.
\end{prop}

\proof. Consider two measurable partitions of unity $D,\,E$ on~$X$. For any $g\in D$ and $h\in E$
define the function $u_{gh} =g\,h\circ\alpha^k$. Let us define  the new partition of unity
$$
C=\{\,u_{gh} \mid g\in D,\ h\in E\,\}.
$$
Set
$$
D_\mu =\{\,g\in D\mid \mu(g)>0\,\}\quad\textrm{and}\quad E_\mu =\{\,h\in E \mid \mu(h)>0\,\}.
$$
If there exists an $h\in E_\mu$ such that $\int_X h\,dm =0$, then $\tau_n(\mu,E) =\tau_{n+k}(\mu,C)
=0$ and so the proposition is proved. In the opposite  case there exist nonnegative functions $f\in
L^1(X,m)$ satisfying the condition $\int_X h\, f\circ\alpha^n\,dm >0$ for all $h\in E_\mu$. For
these $f$, by using the $\alpha$-invariance of $\mu$ and the concavity of the logarithm we obtain
the relations
\begin{align*}
\sum_{u_{gh}\in C}&\mu(u_{gh})\ln\frac{\int_X{u_{gh}}\, f\circ\alpha^{n+k}\,dm}{\mu(u_{gh})}
-\sum_{h\in E_\mu}\mu(h)
\ln\frac{\int_X h\, f\circ\alpha^n\,dm}{\mu(h)} \\[3pt]
&=\sum_{g\in D_\mu}\mu(g)\sum_{h\in E_\mu}\frac{\mu(u_{gh})}{\mu(g)}\ln
\frac{\int_X{u_{gh}}\,f\circ\alpha^{n+k}\,dm\cdot\mu(h)}{
\mu(u_{gh})\cdot\int_Xh\,f\circ\alpha^n\,dm} \\[3pt]
&\le \sum_{g\in D_\mu}\mu(g)\ln\frac{1}{\mu(g)}\sum_{h\in E_\mu}
\frac{\int_X{u_{gh}}\,f\circ\alpha^{n+k}\,dm\cdot\mu(h)}{
\int_Xh\, f\circ\alpha^n\,dm} \\[3pt]
& =\sum_{g\in D_\mu}\mu(g)\ln\frac{\int_Xg\,v\circ \alpha^k\,dm}{\mu(g)},  \qquad \text{where}\quad
v=\sum_{h\in E_\mu}\frac{\mu(h)\,h\,f\circ\alpha^n}{\int_Xh\, f\circ\alpha^n\,dm}.
\end{align*}
Evidently, $\norm v =1$. Varying $f$ we can deduce from here that $\tau_{n+k}(\mu,C)\le
\tau_n(\mu,E) +\tau_k(\mu,D)$. Passing to the infimum over $E$ and $D$ we obtain the desired
inequality from Proposition \ref{..11}. \qed

\section{Proof of the Entropy Statistic Theorem}\label{4..}

Let us fix a natural number $n$ and a measurable partition of unity $D$ on $X$. As before, suppose
that $D_m =\{\,g\in D \mid \int_X g\,dm>0\,\}$ and that $M(D)$ and $M(D_m)$ denote the
finite-dimen\-sional simplexes consisting of all probability measures on $D$ and on $D_m$,
respectively. In this case the simplex $M(D_m)$ is naturally embedded into $M(D)$. For any measure
$\mu\in M(D_m)$ choose a sequence of nonnegative functions  $f_i\in L^1(X,m)$ with unit norms in
such a way that the supremum in \eqref{,,4} is attained at  $\{f_i\}$ and at the same time for all
$g\in D_m$ there exist limits
\begin{equation}\label{,,14}
\mu'(g) =\lim_{i\to\infty}\int_X g\,f_i\circ\alpha^n\,dm.
\end{equation}
Then
\begin{equation}\label{,,15}
\tau_n(\mu,D) =\sup_{\norm f=1}\sum_{g\in D_m}\mu(g)\ln \frac{\int_X
g\abs{f\circ\alpha^n}dm}{\mu(g)} = \sum_{g\in D_m}\mu(g)\ln\frac{\mu'(g)}{\mu(g)}.
\end{equation}
Generally speaking, the correspondence $\mu\mapsto\mu'$ constructed above is many-valued. But for
the convenience of the further representation we now fix some single-valued branch of this
correspondence (i.\,e., to each $\mu\in M(D_m)$ we assign  a unique measure $\mu'$ on $D_m$ of the
form~\eqref{,,14}).

\begin{lemma}\label{..12}
Each summand\/ $\mu(g)\ln\bigl(\mu'(g)\big/\mu(g)\bigr)$ in the right-hand side of\/ \eqref{,,15}
is a bounded function on\/ $M(D_m)$ and it tends to zero as\/ $\mu(g) \to 0$.
\end{lemma}

\proof. It follows from  \eqref{,,14}  that $\mu'(g)\le\norm{A^n}$. So the term
$\mu(g)\ln\bigl(\mu'(g)\big/\mu(g)\bigr)$ is bounded above and has a nonpositive upper limit as
$\mu(g)\to 0$. By Proposition \ref{..6} the function $\tau_n(\mu,D)$ depends continuously on
$\mu\in M(D_m)$. Hence it is bounded on $M(D_m)$. And since all the summands in the right-hand side
of \eqref{,,15} are bounded from above it follows that they are bounded from below as well.

Suppose that $h\in D$ and that $\mu(h) \ln\bigl(\mu'(h) \big/ \mu(h)\bigr)$ does not tend to zero
as $\mu(h)\to 0$. Then there exist a sequence $\mu_i\in M(D_m)$ and an $\eps>0$ such that
$\mu_i(h)\to 0$ and at the same time $\mu_i(h) \ln\bigl(\mu'_i(h) \big/ \mu_i(h)\bigr) <-\eps$.
Replacing the sequence $\mu_i$ by a subsequence we can provide the simultaneous existence of the
limits $\mu(g) =\lim\mu_i(g)$ and $\mu^*(g) =\lim\mu'_i(g)$ as $i\to\infty$. Then
\begin{equation*}
\limsup_{i\to\infty} \tau_n(\mu_i,D) = \limsup_{i\to\infty} \sum_{g\in D_m}
\mu_i(g)\ln\frac{\mu'_i(g)}{\mu_i(g)} \le \sum_{\mu(g)>0}\mu(g)\ln\frac{\mu^*(g)}{\mu(g)}-\eps.
\end{equation*}
It is easy to see that the right-hand side of the latter  inequality does not exceed $\tau_n(\mu,D)
-\eps$, which contradicts the continuity of the restriction of~$\tau_n(\,\cdot\,,D)$ to~$M(D_m)$.
\qed \proofskip

Denote by $\Int M(D_m)$ the set of measures $\mu\in M(D_m)$ that are strictly positive at each
element  $g\in D_m$. For a pair of measures $\mu\in \Int M(D_m)$ and $\nu\in M(D)$ put
\begin{equation}\label{,,16}
\tau_n(\nu,\mu,D) =\sum_{g\in D_m}\nu(g)\ln \frac{\mu'(g)}{\mu(g)},
\end{equation}
where $\mu'(g)$ is defined by~\eqref{,,14}.

\begin{lemma}\label{..13}
For any measure\/ $\mu_0\in M(D_m)$ and any\/ $t>\tau_n(\mu_0,D)$ there exist a neighborhood\/
$O(\mu_0)$ in\/ $M(D)$ and a measure\/ $\mu\in O(\mu_0)\cap\Int M(D_m)$ such that for all\/ $\nu\in
O(\mu_0)$ we have\/ $\tau_n(\nu,\mu,D)<t$.
\end{lemma}

\proof. Let $D_m=\{g_1,\dots,g_k\}$ and $\mu_1$ be the center of the simplex $M(D_m)$. The latter
is defined by $\mu_1(g_i) =1/k$, where $g_i\in D_m$. Consider the family $\mu_\theta =(1-\theta)
\mu_0 +\theta\mu_1$. It lies in $\Int M(D_m)$ provided $\theta\in (0,1]$. Denote by
$O_\theta(\mu_0)$ the set of measures $\nu\in M(D)$ satisfying the inequalities $\abs{\nu(g_i)
-\mu_0(g_i)} <\theta$ for all $g_i\in D_m$. Obviously, $\mu_\theta\in O_\theta(\mu_0)$.

By virtue of Proposition \ref{..6} the function $\tau_n(\mu_\theta,D)$ depends continuously on the
parameter $\theta\ge 0$ and so it is close to $\tau_n(\mu_0,D)$ if $\theta$ is small. So it is
sufficient to prove that if $\nu\in O_\theta(\mu_0)$ then the difference
\begin{equation}\label{,,17}
\tau_n(\nu,\mu_\theta,D) -\tau_n(\mu_\theta,D) = \sum_{g\in D_m} \left(\frac{\nu(g)}{\mu_\theta(g)}
-1\right) \times \mu_\theta(g)\ln\frac{\mu'_\theta(g)}{\mu_\theta(g)}
\end{equation}
uniformly converges to zero as $\theta\to 0$. This can be easily deduced from Lemma~\ref{..12}.
Indeed, if $\mu_0(g)>0$ and $\theta$ is small then the first multiplier of the corresponding
summand in \eqref{,,17} is small while the second multiplier is bounded. And if $\mu_0(g) =0$ then
the first multiplier is bounded while the second one is small. \qed

\begin{lemma}\label{..14}
If\/ $\mu\in\Int M(D_m)$ and the function\/ $f\in L^1(X,m)$ is nonnegative then
\begin{equation*}
\sum_{g\in D_m}\frac{\mu(g)}{\mu'(g)}\int_X g\,f\circ\alpha^n\,dm \le \int_X f\,dm.
\end{equation*}
\end{lemma}

\proof. Without loss of generality we can assume that $\int_X f\,dm =1$. Lemma~\ref{..12} implies
that $\mu'(g)>0$ whenever $\mu(g)>0$. Consider the function
\begin{equation*}
\psi(t) =\sum_{g\in D_m}\mu(g)\ln\frac{(1-t)\mu'(g) +t\int_X g\,f\circ\alpha^n\,dm}{\mu(g)}
\end{equation*}
defined on the segment $[0,1]$. By \eqref{,,14} and \eqref{,,15} it attains  maximum at the point
zero. Hence its derivative at zero is nonpositive:
\begin{equation*}
\psi'(0) =\sum_{g\in D_m}\mu(g)\frac{\int_X g\,f\circ\alpha^n\,dm -\mu'(g)}{\mu'(g)} =\sum_{g\in
D_m}\mu(g)\frac{\int_X g\,f\circ\alpha^n\,dm}{\mu'(g)} -1\le 0.
\end{equation*}
So the lemma is proved. \qed \proofskip

Now we are able to prove the Entropy Statistic Theorem itself. Let us fix a measure $\mu_0\in
M(X,m)$ and an arbitrary number $t>\tau(\mu_0)$. Choose a natural number $n$ and a measurable
partition of unity $D$ on $X$ such that $\tau_n(\mu_0,D) <nt$.

Let  $D_m$ be  the set of all functions $g\in D$ satisfying the inequality $\int_X g\,dm >0$.
Obviously, each $g\in D\setminus D_m$ vanishes $m$-almost everywhere. Hence the function $S_ng(x)
=n\delta_{x,n}(g)$ also vanishes almost everywhere. Thus, $\delta_{x,n}\in M(D_m)$ for almost all
$x\in X$.

Suppose first that $\mu_0\in M(D)\setminus M(D_m)$. Take any neighborhood $O(\mu_0)\subset M(D)$
that has the empty intersection with~$M(D_m)$. In this case $\delta_{x,n}\notin O(\mu_0)$ almost
everywhere. Hence $m\bigl(X_n(O(\mu_0))\bigr) =0$ and so for the situation considered the theorem
is proved.

It remains to consider the case $\mu_0\in M(D_m)$. Then by Lemma~\ref{..13} there exist a
neighborhood $O(\mu_0)$ in $M(D)$ and a measure $\mu\in O(\mu_0)\cap \Int M(D_m)$ such that for all
$\nu\in O(\mu_0)$ the estimate $\tau_n(\nu,\mu,D) <nt$ holds.

Set
\begin{equation}\label{,,18}
\psi(x) =\sum_{g\in D_m}g(x)\ln\frac{\mu(g)}{\mu'(g)},
\end{equation}

\medskip\noindent
where $\mu'(g)$ is defined by \eqref{,,14} and \eqref{,,15}. If we compare \eqref{,,18} with
\eqref{,,16}, it becomes clear that for any natural~$N$ the following equality holds
\begin{equation}\label{,,19}
S_N\psi(x) =-N\tau_n(\delta_{x,N},\mu,D).
\end{equation}

Recall that $\sum_{g\in D_m}g = 1$ almost everywhere. Therefore the convexity of the exponent and
Lemma~\ref{..14} imply that for any nonnegative function $f\in L^1(X,m)$
\begin{equation}\label{,,20}
\int_X e^{\psi}f\circ\alpha^n\,dm\le \int_X\sum_{g\in
D_m}g\,\frac{\mu(g)}{\mu'(g)}\,f\circ\alpha^n\, dm \le \int_X f\,dm.
\end{equation}

Let us introduce the notation
$$
\psi_k =\psi+\psi\circ\alpha^n+\dots+\psi \circ\alpha^{n(k-1)}.
$$
Then
$$
e^{\psi_k} f\circ\alpha^{nk} =e^{\psi} \bigl(e^{\psi_{k-1}}f\circ \alpha^{n(k-1)}
\bigr)\circ\alpha^n.
$$

\medskip\noindent
Applying $k$ times estimates \eqref{,,20} to the latter equality  we obtain
\begin{equation}\label{,,21}
\int_Xe^{\psi_k}f\circ\alpha^{nk}\,dm\le \int_Xf\,dm.
\end{equation}

By construction, $\psi(x)$ is essentially bounded and thus there exists a constant $C$ such that
$\abs{\psi(x)} \le C$ for almost all $x\in X$. Consider a natural number $N>n$. For this number
choose  an integer $k$ such that $N\in [n(k+1),\,n(k+2)]$. Then
\begin{equation}\label{,,22}
S_N\psi\le S_{nk}\psi+2nC =\sum_{i=0}^{n-1}\psi_k\circ\alpha^i +2nC.
\end{equation}

If $x\in X_N(O(\mu_0))$ then $\delta_{x,N}\in O(\mu_0)$ and by the choice of $O(\mu_0)$ we have the
estimate $\tau_n (\delta_{x,N},\mu,D) <nt$. In this case \eqref{,,19} and \eqref{,,22} imply
\begin{equation*}
Nt>\frac{N}{n}\tau_n(\delta_{x,N},\mu,D) = -\frac{1}{n}S_N\psi(x)\ge -\frac{1}{n}\sum_{i=0}^{n-1}
\psi_k\circ\alpha^i(x) -2C
\end{equation*}
for almost all $x\in X_N(O(\mu_0))$. Whence,
\begin{equation*}
\int_{X_N(O(\mu_0))} f\circ\alpha^N\,dm \le
e^{Nt+2C}\int_{X}\exp\biggl\{\frac{1}{n}\sum_{i=0}^{n-1}
\psi_k\circ\alpha^i\biggr\}f\circ\alpha^N\,dm.
\end{equation*}
By virtue of H\"older's inequality the latter  integral does not exceed
\begin{equation*}
\prod_{i=0}^{n-1}\left(\int_X e^{\psi_k\circ\alpha^i} f\circ\alpha^N\,dm\right)^{\!1/n}\le
\prod_{i=0}^{n-1} \left(\norm{A}^i\!\int_X e^{\psi_k}f\circ\alpha^{N-i}\,dm \right)^{\!1/n},
\end{equation*}
and by \eqref{,,21} one has
\begin{equation*}
\norm{A}^i\!\int_X e^{\psi_k}f\circ\alpha^{N-i}\,dm \le \norm{A}^i\!\int_X
f\circ\alpha^{N-i-nk}\,dm \le \norm{A}^{N-nk}\!\int_X f\,dm.
\end{equation*}
Now the Entropy Statistic Theorem follows from the final  three estimates. \qed

\end{document}